\documentclass{amsart}

\usepackage{amssymb}
\usepackage{amscd}

\title{Galois points for the Dickson--Guralnick--Zieve curve}
\author{Satoru Fukasawa}

\subjclass[2010]{14H50, 14H37, 14G05}
\keywords{Galois point, plane curve, rational point}
\address{Department of Mathematical Sciences, Faculty of Science, Yamagata University, Kojirakawa-machi 1-4-12, Yamagata 990-8560, Japan}
\email{s.fukasawa@sci.kj.yamagata-u.ac.jp} 
\thanks{The author was partially supported by JSPS KAKENHI Grant Number 16K05088.}  

\newtheorem{theorem}{Theorem}

\newtheorem{fact}{Fact}
\newtheorem{problem}{Problem} 

\theoremstyle{definition}

\begin{document}
\begin{abstract} 
The Dickson--Guralnick--Zieve curve over a finite field has been studied recently by Giulietti, Korchm\'{a}ros and Timpanella in several points of view. 
In this short note, the distribution of Galois points for this curve is determined. 
As a consequence, a problem posed by the present author in the theory of Galois point is modified. 
\end{abstract}
\maketitle

\section{Introduction}  
The Dickson--Guralnick--Zieve curve over a finite field $\mathbb{F}_q$ has a large automorphism group and a positive $p$-rank, and hence, which is an important class in the study of automorphism groups of algebraic curves. 
The several good properties of this curve has been studied recently by Giulietti, Korchm\'{a}ros and Timpanella (\cite{gkt}). 
In this short note, we consider Galois points for this curve over $\overline{\mathbb{F}}_q$ (see \cite{miura-yoshihara, yoshihara} for the definition of Galois point). 
The set of all Galois points for a plane curve $\mathcal{C}$ on the projective plane is denoted by $\Delta(\mathcal{C})$. 
It would be good to obtain a new example of a plane curve $\mathcal{C}$ with large $\Delta(\mathcal{C})$. 

Our main result is the following. 

\begin{theorem}
For the Dickson--Guralnick--Zieve curve $\mathcal{C} \subset \mathbb{P}^2$, 
$$ \Delta(\mathcal{C})=\mathbb{P}^2(\mathbb{F}_q). $$
In particular, the number of Galois points for $\mathcal{C}$ is exactly $q^2+q+1$. 
\end{theorem}

Since this result gives a negative answer to the problem \cite[Problem 1]{fukasawa2} posed by the present author, this is modified as follows. 

\begin{problem} 
Let $\mathcal{C}$ be a plane curve over $\mathbb F_q$. 
Assume that $\Delta(\mathcal{C})=\mathbb P^2(\mathbb F_q)$. 
Then, is it true that $\mathcal{C}$ is projectively equivalent to the Hermitian, Ballico--Hefez or Dickson--Guralnick--Zieve curve?
\end{problem}

\section{Proof}
Let $\mathbb F_q$ be a finite field with $q \ge 2$. 
The Dickson--Guralnick--Zieve curve $\mathcal{C} \subset \mathbb{P}^2$ is the plane curve defined by $F(x, y, z)=D_1(x, y, z)/D_2(x, y, z)$, where 
$$ D_1=\left|\begin{array}{ccc}
x & x^q & x^{q^3} \\
y & y^q & y^{q^3} \\ 
z & z^q & z^{q^3} 
\end{array} \right| \  \mbox{ and } \ 
D_2=\left|\begin{array}{ccc}
x & x^q & x^{q^2} \\
y & y^q & y^{q^2} \\ 
z & z^q & z^{q^2} 
\end{array} \right|. $$ 
According to \cite[Lemma 4.2 and Proposition 4.7]{gkt}, $F$ is a homogeneous polynomial of degree $q^3-q^2$ over $\mathbb{F}_q$ and is irreducible over the algebraic closure $\overline{\mathbb F}_q$.

It is remarkable that the projective linear group $PGL(3, \mathbb{F}_q)$ acts on $\mathcal{C}$ (\cite[Lemma 4.1]{gkt}). 
Therefore, the matrices
$$\sigma_{\gamma, \beta}=\left(\begin{array}{ccc}
1 & 0 & 0 \\
\gamma & 1 & \beta \\ 
0 & 0 & 1 
\end{array} \right) \ \mbox{ and } \ 
\tau_{\mu}=\left(\begin{array}{ccc}
1 & 0 & 0 \\
0 & \mu & 0 \\ 
0 & 0 & 1 
\end{array} \right) \in PGL(3, \mathbb{F}_q)$$
act on $\mathcal{C}$, where $\gamma, \beta \in \mathbb{F}_q$ and $\mu^{q-1}=1$. 
Note that a rational function $x/z$ is fixed by the actions of $\sigma_{\gamma, \beta}$ and $\tau_{\mu}$. 
This implies that $\pi_{P} \circ \sigma_{\gamma, \beta}=\pi_{P}$ and $\pi_{P} \circ \tau_{\mu}=\pi_{P}$, where $P=(0:1:0)$ and $\pi_P$ is the projection from $P$. 
It follows that $P$ is a Galois point not contained in $\mathcal{C}$. 
Considering the action of $PGL(3, \mathbb{F}_q)$, we infer that $\mathbb{P}^2(\mathbb{F}_q) \subset \Delta(\mathcal{C})$ holds. 

We would like to show that $\Delta(\mathcal{C}) \subset \mathbb{P}^2(\mathbb{F}_q)$ holds. 
When $q=2$, the curve $\mathcal{C}$ is given by 
$$ F(x, y, z)=(x^2+xz)^2+(x^2+xz)(y^2+yz)+(y^2+yz)^2+z^4$$
(\cite[Remark 1]{gkt}). 
In this case, it is known that the claim follows (\cite[Theorem 4]{fukasawa1}). 
Hereafter, we assume that $q>2$. 
The following assertions have been obtained by Giulietti, Korchm\'{a}ros and Timpanella \cite[Lemmas 4.5, 4.6 and 8.4 and Corollary 4.8]{gkt}. 

\begin{fact} \label{properties} 
Let $q>2$, let $\mathcal{C}$ be the Dickson--Guralnick--Zieve curve, and let $r: \hat{\mathcal{C}} \rightarrow \mathcal{C}$ be the normalization. 
\begin{itemize}
\item[(1)] ${\rm Sing}(\mathcal{C})=\mathbb{P}^2(\mathbb{F}_{q^2}) \setminus \mathbb{P}^2(\mathbb{F}_q)$. 
\item[(2)] $r$ is bijective. 
\item[(3)] For any point $Q \in {\rm Sing}(\mathcal{C})$ and any line $\ell \ni Q$, ${\rm ord}_Q\ell=q-1$ or $q$. 
Furthermore, if ${\rm ord}_Q\ell=q$, then $\ell$ is defined over $\mathbb{F}_q$. 
\item[(4)] For any point $Q \in \mathcal{C} \setminus {\rm Sing}(\mathcal{C})$ and the tangent line $T_Q$ at $Q$, ${\rm ord}_QT_Q \ge q$. 
\end{itemize}
\end{fact}

Let $P \in \Delta(C)$. 
Note that any point $Q \in {\rm Sing}(\mathcal{C})$ is a ramification point of the projection $\pi_P$, and $e_Q=q-1$ or $q$, by Fact \ref{properties}(2)(3).
If $e_Q=q-1$, then it follows from Fact \ref{properties} (3)(4) and \cite[III.7.2]{stichtenoth} that the line $\overline{PQ}$ passing through $P$ and $Q$ intersects with $\mathcal{C}$ at two or more $\mathbb{F}_{q^2}$-points and hence, the line $\overline{PQ}$ is $\mathbb{F}_{q^2}$-rational.  
If $e_Q=q$, then the line $\overline{PQ}$ is $\mathbb{F}_{q}$-rational, by Fact \ref{properties}(3). 
Since $P$ is the intersection point given by some two $\mathbb{F}_{q^2}$-lines, it follows that $P \in \mathbb{P}^2(\mathbb{F}_{q^2})$. 

Assume that $P \in \mathbb{P}^2(\mathbb{F}_{q^2}) \setminus \mathbb{P}^2(\mathbb{F}_q)={\rm Sing}(\mathcal{C})$.  
We consider the tangent line $T_P$ at $P$, i.e. ${\rm ord}_PT_P=q$. 
Then, $T_P$ is $\mathbb{F}_q$-rational and contains a Galois point in $\mathbb{P}^2(\mathbb{F}_q)$ as above. 
It follows from \cite[III.7.2]{stichtenoth} that there exists a point $Q \in \mathcal{C} \cap T_P$ with $Q \ne P$ such that ${\rm ord}_QT_P=q$. 
Then, $e_Q=q$ for the projection $\pi_P$. 
It follows from \cite[III.7.2]{stichtenoth} that $e_Q=q$ divides $\deg \pi_P=q^3-q^2-q+1$. 
This is a contradiction. 
The assertion $P \in \mathbb{P}^2(\mathbb{F}_q)$ follows.

\end{document}